\documentclass[11pt,a4]{amsart}

\usepackage{amssymb,epsfig}
\usepackage{amsmath}
\newif\ifpdf
\ifx\pdfoutput\undefined
\pdffalse 
\else
\pdfoutput=1 
\pdftrue
\fi

\ifpdf
\usepackage[pdftex]{graphicx}
\else
\usepackage{graphicx}
\fi

\newtheorem{theorem}{Theorem}
\newtheorem{proposition}{Proposition}
\newtheorem{definition}{Definition}
\newtheorem{remark}{Remark}

\def\rit{\mathbb{R}} 
\def\R{\mathbb {R}}

\def\nit{\mathbb{N}} 
\def\cit{\mathbb{C}} 
\def\C{\mathbb {C}}
\def\zit{\mathbb{Z}}

\def\Z{{\mathbb Z}}

\def\CD{\mathcal{D}}
\def\CI{\mathcal{I}}
\def\CL{\mathcal{L}}
\def\CT{\mathcal{T}}

\newcommand{\supe}{\mathop{\sup}}

\newcommand{\Sup}[1]{\supe\limits_{#1}}

\newcommand{\ds}{\displaystyle}

\def \bx {\hspace{2.5mm}\rule{2.5mm}{2.5mm}}
\def \cqfd {\bx}
\def \O {\mathcal{O}}

\def \Sp {{\rm Sp}}

\def \Entier{{\rm \lfloor}} 
\def \Im {{\rm Im}}

\def \ln {{\rm ln}}
\def \im {{\rm Im\,}} 
\def \re {{\rm Re\,}}

\def \a {\alpha}

\def \po {\rm{o}}

\def \det {{\rm det}}

\title{ Perturbations of\\
 non self-adjoint  Sturm-Liouville problems, with  applications to  harmonic oscillators}

\author{Laurence Nedelec$\,^\star$} 
\address{L. Nedelec,  L.A.G.A., Institut Gali\-l\'ee, Univer\-sit\'e de Paris Nord, av. J.B. Cle\-ment,
F-93430 Villetaneuse, France\\
 IUFM de l'academie de Rouen, France}
\email{ nedelec@math.univ-paris13.fr}
\thanks{$\,^\star$  Supported in part by the FIM of ETHZ}

\begin{document}


\begin{abstract} 
We study the behavior of the limit of the spectrum of a non self-adjoint
 Sturm-Liouville operator with  analytic potential as 
 the semi-classical parameter $h\to 0$. We get a good description of the spectrum and  limit spectrum near $\infty$. We also study the action of one  special perturbation of the operator (adding a Heaviside function), and prove  
that the limit spectrum is very unstable. 
 As an illustration we describe the limit spectrum as  $h\to 0$ for $P^h=-h^2\Delta+i x^2$ and the effect of this perturbation. 
\end{abstract}
\maketitle

%
\section{Introduction}

This paper is devoted to  non self-adjoint Sturm-Liouville problems. We study  the spectrum  of the   1-dimensional, semiclassical Schr\"odinger operator on $L^2([-1,1])$ with Dirichlet boundary condition, given by
\begin{equation}
H^h=-h^2 \frac d {dx^2} +V(x).
\label{Hi}
\end{equation}
The potential $V$ is a complex valued function  on $[-1,1]$, which extends holomorphically to some domain in $\C$.
The boundary value at $\pm 1$ play no special rule but are fixed to avoid more notation. 

The study of such operators is motivated  by  the Orr-Som\-merfeld equation
with linear profile \cite{sh} or by the non linear Zakharov-Shabat eigenvalue problem, cf.  work of Miller \cite{mi}. 

As an application  we will focus on the case where 
\begin{equation}
V(x)=i x^2,
\label{Vquad}
\end{equation}
and shall write $P^h$ for the corresponding operator. From this one could also
study the slightly  more general case 
$-h^2 \frac d {dx^2} +e^d x^2$, $d\in\cit$ using a change of variable.
The spectrum of this operator on $\rit$  was analyzed by Davies \cite{davies},
cf. also the recent work of Hitrik \cite{hit}.

It is well known that the spectrum of a non self-adjoint operator is unstable
under perturbation of the operator. This motivates the introduction of the
pseudo-spectrum, which has now been studied by many people, particularly
Trefethen (who maintains the web archive http:web.comlab.ox.ac.uk/projects/pseudospectra)
and Davies \cite{dav},\cite{dave}; we note also the recent paper of 
Denker, Sj\"ostrand and  Zworski \cite{speudo}. 

We  also consider the following perturbation of $H$: 

For $\beta\in (-1,1)$ and $\delta\geq 0$, let $H_{\delta,\beta}$
\begin{equation}
H^h_{\delta,\beta}=-h^2 \frac d {dx^2} +V_{\delta,\beta}(x), \quad V_{\delta,\beta}(x)=
\left\{
\begin{array}{l}
V(x)+i\delta,\quad x>\beta\\
V(x)-i\delta,\quad x<\beta
\end{array}
\right., 
\label{Hdelta}
\end{equation}
with domain 
\begin{equation}
\CD_{\delta,\beta}=\{u\in L^2([-1,1]), u''\in L^2([-1,1]), u(-1)=u(1)=0, \}.
\label{domain}
\end{equation}
Notice that $H^h_{0,\beta}=H^h$.

Redparth \cite{red} has obtained similar results for arbitrary piecewise
linear complex-valued potentials, a class which includes perturbations as above 
of the operator $-h^2\Delta+ix$.

We  focus on the semi classical limit, and denote by  
$\lim_h\Sp(H^h_\delta)$ the set of all values of $E$ which arise as limits  as $h\to 0$ of some sequence of eigenvalues  of $ H^h_\delta$.
 
We give  a geometric criterion for $E$ not to  belong to  $\lim_h\Sp(H^h_\delta)$:   namely 
$\lim_h\Sp(H^h_\delta)\subset \CT_\delta^c,\quad \CT_\delta^c=\cit\subset\CT$ where $\CT_{\delta}$ is the set of $E\in\C$ such that there exists
a progressive path from  $-1$ to $1$ with respect to $V_{\delta,\beta}$. ( We
define this concept later). This is done in section \ref{semi},  see Theorems \ref{inclus}. 
 This criterion is not sharp, and also not so easy to use. 
It can be proved using either  ellipticity  or  exact  WKB .

As a second step we give an alternate criterion for $E$ to belong to
\par\noindent $\lim_h\Sp(H^h_\delta)$, when $E$ is large enough, which is sharp and
computable, see Theorem \ref{reci} .

As an application we prove that the perturbation $\delta$ changes the spectrum
near  $\infty$ quite drastically, see Theorem \ref{curveper}.
The pseudo spectrum is defined to be the limit as $\delta\to 0$ of the union of spectrum over all pertubations of $H^h$ of size $\delta$, and we show that  this perturbation is sufficient to reach
the entire pseudo spectrum of $H^h$, see section \ref{speudospec}.
We note also the work  in progress of M. Hager that studies the effect of different pertubations on the spectrum. 

We use the geometric criterion for  $P^h$ to  characterize the full set 
$\CT^c$, see Theorem \ref{forme}, using special properties of the potential,
and see that $\CT^c$ forms a $Y$ shape. After this article was completed, we
learned that some of our results and methods are similar to those of Shkalikov
\cite{sh}, \cite{sha}.
The algorithm proposed by Miller in \cite{mi} to study the spectrum is
the same as the one we use here, though he applies it to a different operator.   

We also obtain a description of the large spectrum of $H^h$ which is near to
\par\noindent $\lim_h\Sp(P^h)$ see Theorem \ref{quanti}. 
In the  self-adjoint case; this can be found in the book of Marchenko \cite{ma} or the paper of Kappeler and M\"ohr  \cite{KM}, and in  some special  non self-adjoint cases in work of Carlson, Threadgill and Shubin \cite{ca}. 
However, these results apply for  potentials $V$ with  $V'\notin L^1_{loc}$,
hence do not pertain to  $V_{\delta,\beta}$.

Even though   it is not self-adjoint, the operator $H^h_{\delta,\beta}$ has discrete spectrum. Indeed, $V_{\delta,\beta}$ is  compact relative to $\Delta$, and from Weyl's theorem it follows that  its essential  spectrum is empty. In particular $H^h_{\delta,\beta}$ has no residual spectrum.

%

\section{Generals tools}

In this section we are interested in the Dirichlet  eigenvalues  of the operator $H^h$ on $L^2([-1,1])$ defined by 
\begin{equation}
\left\{
\begin{array}{ll}
H^hu=Eu,\\
u(-1)=u(1)=0.
\end{array}
\right .
\label{schr}
\end{equation}

\begin{definition}\sl A point $x_0$ in the complex plane is a turning point of order $k$ for  the operator $H^h-E$ if $V(x)-E$ vanishes to order $k$ at $x=x_0$.
\label{turningpoint}
\end{definition}

\begin{definition}\sl
Let $x$ be a point in the complex plane. The Stokes line starting from $x$ is the set 
$$
\CL(x)=\{ y\in \C, \re\int_x^y (V(t)-E)^{1/2} \,dt=0\}.
$$
\label{defLS}
\end{definition}
This is independant of the branch of the square root.
Denote by $S_{x,y}(E)=\int_x^y (V(t)-E)^{1/2} \,dt$.

Notice that the Stokes lines are integral curves of the Stokes vector field \begin{equation}
s(x)=i\;\overline{(V(x)-E)^{1/2}}.
\label{stokesfield}
\end{equation}

The local structure of Stokes lines can be easily investigated (see
e.g. \cite{fe}, \cite{fr}).  Away from  turning points or singularities
of $V$,  the Stokes lines form a non singular analytic foliation.
 Indeed, if $\Omega\subset \C$ does not contain any singularities of $(V-E)^{1/2}$,  then $x\mapsto z(x,E)=\int^x_{x_0}(V(t)-E)^{1/2}\,dt$ is an analytic diffeomorphism.
Notice also that from a simple turning point $x_0$ emanates  three Stokes lines,  each making an angle of $2\pi/3$ with any other at $x_0$.

Using ideas from quantum  resonance theory, we shall consider as in \cite{fr}  some distorted operators associated to $H$.
Let $\gamma:[-1,1]\ni t \to \C$ be a smooth simple path in the complex plane $\gamma(-1)=-1$ and $\gamma (1)=1$.
Since $V$ is analytic, we can define an operator $H^{\gamma,h}$ on
$L^2([-1,1])$ as the restriction of  $H^h$  to $\gamma([-1,1])$. One computes that
$$H^{\gamma,h}= -h^2 ( \frac 1 {\gamma'(x)} \frac d  {dx}) ( \frac 1 {\gamma'(x)} \frac d  {dx}) +V(\gamma(x)).    $$

The correspondence between the Stokes geometry and spectral properties of $H^h$ is given in the following result  from \cite{fr}.

\begin{proposition}\sl  
The operator $H^{\gamma,h} -E$ is elliptic if and only if $\gamma$ is transverse to the Stokes lines.
\label{ellipticite}
\end{proposition}

{\sl Proof: } The principal symbol of  $H^{\gamma,h} -E$ is 
$$
h^\gamma(t,\tau)=\frac 1{(\gamma'(t))^2}(\tau^2+(\gamma'(t))^2(V(\gamma(t))-E)).
$$
We set  $(V(\gamma(t))-E)^{\frac 1 2}= r e^{i\theta}$ and $\gamma'(t)=e^{i\theta'}$  (one  can always suppose that $\vert\gamma'(t)\vert=1$).
Then the path $\gamma$ is transversal to the Stokes lines if and only if 
$$
\det (\gamma'(t), s(\gamma(t)))\neq 0,
$$
where $s(x)$ is the Stokes field  defined in (\ref{stokesfield}). This condition is the same as
$$
\re \gamma'(t)(V(\gamma(t)-E)^{1/2}\neq 0, 
$$
or finally
$$
\cos(\theta+\theta')\neq 0.
$$
Therefore, the transversality condition is equivalent to
$\theta+\theta'\not\equiv \pi/2  \mod [\pi]$.
Since $\im ({\gamma'(t)}^2 h^{\gamma})=r^2\sin  2(\theta+\theta')$  and $\re( {\gamma'(t)}^2  h^{\gamma})=\tau^2+r^2\cos 2(\theta+\theta')$,  the transversality condition  is equivalent to the invertibility of $h^\gamma$, as stated in the proposition. \cqfd

\medskip

We deduce the following proposition, but  this will be improved in Theorem
\ref{mainbis} below, where we control the dependence on $E$ and introduce the perturbation.
\begin{proposition}\sl  
For $h$ small enough, $E\notin \Sp(H^h)$  if there exists  $\gamma$  form $-1$ to $1$ which is transverse to the Stokes lines.
\label{inver}
\end{proposition}
{\sl Proof: } 
General properties of  distorted analytic operators imply that  $E\in \Sp(H^h)$ is equivalent to $E\in \Sp(H^{\gamma,h})$. 
Let $\gamma$ be a path  transverse to the Stokes lines. Then by Proposition \ref{ellipticite},
$H^{\gamma,h} -E$ is elliptic and invertible for sufficiently small $h$ (depending on $E$). 
\cqfd

%

\section{The semi-classical limit of the spectrum }
\label{semi}
The section is devoted to the proof of the  Theorem \ref{mainbis} and a direct application Theorem \ref{inclus}.
The next sections \ref{large} and \ref{largeper}  are devoted to the proof of
a converse
of Theorem \ref{inclus}. 

We now introduce  some notation.
Denote by $\lim_h \Sp(H^h)$ the set of all $E_0$ such that 
there exists a sequence $h_j\to 0$  with  $E_{h_j}\in\Sp(H^{h_j})$ and  
$E_{h_j}\to E_0$. For $a\in\cit$,
we denote by $z_{0}(x)$ and $z_{l,r}(x)$ the  action integrals
\begin{equation}
\begin{array}{c}
\ds z_{0}(x)=\int_{\gamma(-1,x)} (V(t)-E)^{1/2}\,dt,\;
z_{l}(x)=\int_{\gamma_{l}(\beta,x)} (V(t)-i\delta-E)^{1/2}\,dt,
\\
\ds z_{r}(x)=\int_{\gamma_{r}(\beta,x)} (V(t)+i\delta-E)^{1/2}\,dt,
\end{array}
\label{actions}
\end{equation}
where $\gamma(-1,x)$ and $\gamma_{l,r}(\beta,x)$ are paths from $-1$ or $\beta$ to
$x$ on the Riemann surface associated to the determination of $\Sigma_{0}$ of $t\mapsto (V(t)-E)^{1/2}$ and
$\Sigma_{l,r}$ associated to the determination of $t\mapsto (V(t)\pm i\delta-E)^{1/2}$, $E\in\cit$ respectively.

We denote by $\CT_{\delta}$ the set of $E\in\C$ such that there exists $\gamma_{l}(-1,\beta)$ on which the function $t\mapsto\re(z_l(\gamma_{l}(t) ))$ is
strictly monotone, and  a path $\gamma_{r}(\beta,1)$ on which the function
$t\mapsto\re(z_r(\gamma_{r}(t) ))$ is strictly monotone. 
We will call such a path 'progressive' .
For $\delta=0$ we take a different definition: $E\in\CT_0$ if and only if  there exists a path $\gamma$ from $-1$  to
$1$ on which the function $t\mapsto \re(z_0(\gamma(t))))$ is
strictly monotone.

Finally, for $E\in \CT_\delta$,  we denote by
\begin{equation}
d(E,\CT_\delta^c,\gamma)=\inf_{l,r}\inf_{x\in \gamma^{l,r}} |\re \partial_x z^{l,r}(x)| \quad \mbox{ if }  \delta\neq 0
\end{equation}
\begin{equation}
d(E,\CT_0^c,\gamma)=\inf_{x\in \gamma} |\re \partial_x z^0(x)| 
\label{dist}
\end{equation}
where $\gamma$  is a progressive path associated to $E$ ($\gamma=\cup_{l,r}\gamma^{l,r}$ if $\delta\neq 0$).

\begin{remark}\sl 
\label{remclose}
Notice that $\CT_\delta$ is open. More precisely, if  $c$ is small enough there exists $C>0$,  such that if $E_0$ belongs to $\CT_\delta$  then, for $E$ such that 
$|E-E_0|\leq c d(E_0,\CT_\delta^c,\gamma)$, we get $E\in \CT_\delta$   and $d(E,\CT_\delta^c,\gamma)\geq C d(E_0,\CT_\delta^c,\gamma)$.
This follows from the relation 
$$
\sqrt{V_{\delta,\beta}-E}=\sqrt{V_{\delta,\beta}-E_0}+O(\frac {E-E_0}{\sqrt{V_{\delta,\beta}-E_0}}),
$$
which shows that a progressive path for $E_0$  is also  a progressive path for $E$.
 \end{remark}

\begin{theorem}
\label{mainbis}\sl 
Let $\delta\geq 0$ be  small enough and  $E\in \CT_{\delta}$. Then  $E\notin
  \Sp(H^h_\delta)$ as soon as $h\leq  d(E,\CT^c_{\delta},\gamma)^7$.
  (The reason for the exponent $7$ will emerge at the end of the proof.)
\end{theorem}

From Remark \ref{remclose} and Theorem \ref{mainbis} we obtain the 
\begin{theorem}\label{inclus}\sl 
If $E_0\in\lim_h\Sp(H^h_\delta)$, then  $E_0\in \CT_\delta^c$. 
\end{theorem}

{\sl Proof:} 
Suppose that $E_0\in\CT_\delta$ and $\gamma$ is a progressive path associated  to $E_0$. Let $E_h$ be a sequence in the spectrum of $H^h_\delta$ tending to $E_0$. By the remark above, $E_h$ belongs to $\CT_\delta$ as soon as $h$ is small enough, and we get $d(E_h,\CT_\delta^c,\gamma)\geq \tilde C_0$.
Theorem \ref{mainbis} shows that $E_h$ is not in the spectrum  of $H^h_\delta$ as soon as $h$ is small. This contradicts the hypothesis, so $E_0\in \CT_\delta^c$.
\cqfd

\medskip
{\sl Proof of Theorem \ref{mainbis}: } 
We denote by $(u^l_+,u^l_-)$ a basis of the space of  solutions of the equation
\begin{equation}
(H^h_\delta-E)u=0,
\label{eqvp}
\end{equation}
in $L^2([-1,\beta ])$, and  by $(u^r_+, u^r_-)$  a basis  of the space of solutions of (\ref{eqvp}) in $L^2([\beta,1])$.
A function $f\in L^2([-1,1])$ is an eigenfunction of $H_\delta$ with eigenvalue $E\in \C$, if and only if   
$f(1)=f(-1)=0$, and both $f$ and $\partial_x f$ are continuous at $x=\beta $. Writing 
\begin{equation}
f=\alpha^l_+u^l_++\alpha^l_-u^l_-\mbox{ on } [-1,\beta] \quad f=\alpha^r_+u^r_++\alpha^r_-u^r_-\mbox{ on } [\beta,1] ,
\label{solglobale}
\end{equation}
Then $f$ is an eigenfunction of $H^h_\delta$ with eigenvalue $E\in \C$ if and only if
\begin{gather}
\begin{split}
\left(
\begin{array}{cccc} 
u^l_-(-1)&u^l_+(-1)&0&0\\
 u^l_-(\beta)&u^l_+(\beta)&-u^r_-(\beta)&-u^r_+(\beta)\\
\partial_x{u^l}_-(\beta)&\partial_x{u^l}_+(\beta)&-\partial_x{u^r}_-(\beta)&-\partial_x{u^r}_+(\beta)\\
0&0&u^r_-(1)&u^r_+(1)\end{array}\right)
\left(\begin{array}{l}\alpha^l_+\\\alpha^l_-\\\alpha^r_+\\
\alpha^r_-\end{array}\right)=0.
\end{split}
\label{defI}
\end{gather}

Therefore $E$ belongs to $\Sp (H^h_\delta)$, the spectrum of $H^h_{\delta}$,  if and only if
\begin{equation}
\det\, \CI(\beta,E)=0,
\label{detI}
\end{equation}
where $\CI(\beta,E)$ is the matrix appearing here.
\begin{remark}
For $\delta=0$, the proof is slightly different, the corresponding matrix is 
\begin{gather}
\begin{split}
\CI(\beta,E)=\left(
\begin{array}{cccc} 
u_-(-1)&u_+(-1)\\
 u_-(1)&u_+(1)\end{array}\right),
\end{split}
\label{defIdelta}
\end{gather}
We leave details to the reader and in the following, treat only the case
$\delta\neq 0$.
\end{remark}

 In the sequel, we shall compute this determinant for two particular bases $(u^l_+,u^l_-)$ and $(u^r_+, u^r_-)$, namely for complex WKB solutions  as in  \cite{gg}, or \cite{fr}.  
 Suppose $E\in T_{\delta}$. Then there is a path
 $\gamma_l=\gamma_{l}(-1,\beta)$ transversal to the Stokes lines $\CL(y)$,
 and we can suppose that $\gamma_{l}(-1,\beta)$ is of type $+$, that is
 $t\mapsto z_{l}(\gamma_{l}(t))$ is strictly increasing. Indeed if it is not
 the case, the path on the other sheet of the Riemann surface $\Sigma_{l}$
 associated to the definition of the square root with the same projection on $\C$ as 
$\gamma_{l}(-1,\beta)$ is then of type $+$. (We have $z_{l}(-1)<0$.)
We define two independent complex WKB solutions $w_{\pm}$ of the equation 
$(H^h_\delta-E)u=0$ on the interval $[-1,\beta]$.These have the form
\begin{gather}
\begin{split}
w^l_\pm : x\mapsto (V(x)-i\delta-E)^{-\frac 1 4} e^{\pm z_l(x)/h}W^l_\pm(z_l(x))
\\
=({\partial_x z_l (x)})^{-\frac 1 2}(x)e^{ \pm z_l(x)/h} W^l_\pm (z_l(x)).
\label{wkbsoll}
\end{split}
\end{gather}

Similarly, there exists a path $\gamma_r=\gamma_r(\beta,1)$ transverse to the
Stokes lines of type $+$ , (We have $z_r(1)>0$ and we have again  two independent complex  WKB solutions  of the equation 
$H_\delta w_\pm=0$ on the interval $[\beta,1]$ :
\begin{gather}
\begin{split}
w^r_\pm : x\mapsto (V(x)+i\delta-E)^{-\frac 1 4} e^{\pm z_r(x)/h}W^l_\pm(z_r(x))\\
=({\partial_x z_r}(x))^{-\frac 1 2}e^{ \pm z_r(x)/h} W^r_\pm (z_r(x)).
\label{wkbsolr}
\end{split}
\end{gather}
The functions $W^{l,r}_\pm$ are convergent series  of the form
\begin{equation}
W^{l,r}_\pm=\sum_{n=0}^\infty W_{n,\pm}^{l,r},
\label{wlrn}
\end{equation}
where the $W_{n,\pm}^{l,r}$ are determined  by the recurrence relations:
\begin{equation}
\left\{
\begin{array}{l}
(\partial_z \pm\frac 2h)W^{l,r}_{2n+1,\pm}=-H^{l,r}W^{l,r}_{2n,\pm}\\[8pt]
\partial_z W^{l,r}_{2n,\pm}=-H^{l,r}W^{l,r}_{2n-1,\pm}
\end{array}
\right .
\end{equation}
with initial data
\begin{equation}
\left\{
\begin{array}{ll}
W^{l,r}_{0,\pm}=1,&\\[8pt]
W^{l}_{n,+}(z_l(-1))=0,\ W^{r}_{n,-}(z_r(1))=0, & n\geq 1,\\[8pt]
W^{l}_{n,-}(z_l(\beta))=0,\ W^{r}_{n,+}(z_r(\beta))=0, & n\geq 1,
\end{array}
\right .
\end{equation}
We have set here
\begin{equation}
 H^{l,r}(z_{l,r}(x))=-\frac {1} 4 \frac {\partial_xV(x)}{(\partial_x z_{l,r}(x))^3}\cdotp
\label{H}
\end{equation}
 These equations can be written in integral form as
\begin{equation}
\left\{
\begin{array}{l}
W_{2n+1,\pm}^{r,l}= I_\pm^{l,r}(W_{2n,\pm}^{r,l}),\\[8pt]
 W_{2n+2,\pm}^{r,l}= J^{l,r}(W_{2n+1,\pm}^{r,l})
 \end{array}
\right .
 \end{equation}
where
\begin{equation}
\left\{
\begin{array}{l}
\ds
I_\pm^{l,r}(v)(z)=-\int_{\tilde\gamma_{l,r}(z)}e^{\pm2(u-z)/ h} H^{l,r}(u)v(u)\,du
\\[12pt]
\ds
J^{l,r}(v)(z)=-\int_{\tilde\gamma_{l,r}(z)} H^{l,r}(u)v(u)\,du 
\end{array}
\right .
\label{IetJ}
\end{equation}
and $\tilde\gamma_{l,r}(z)$ is the image by $x\mapsto z_{l,r}(x)$ of
$\gamma_{l,r}$ from the initial point  $-1$ to $z$ for $l$ and  from $\beta$ to $z$ for $r$. 

Now we have the following estimates
\begin{gather}\begin{split}
\Sup{z\in\tilde\gamma_{l,r}}|I^{l,r}_+(v)(z)| 
&\leq  \Sup{z\in\tilde\gamma_{l,r}} |v(z)| 
\Sup{z\in\tilde\gamma_{l,r},z=z^{l,r}(x)}
\big\{ |H^{l,r}(z)|
\frac {|\partial_x z_{l,r}(x)|}{|\re \partial_x z_{l,r}(x)| } \big\} \times
\\& |\int_{\tilde\gamma_{l,r}}e^{\re (\frac {2(u-z)} h)}\re(\,dz)|\\&
\leq  h \Sup{z\in\tilde\gamma_{l,r}} |v(z)|
\Sup{z\in\tilde\gamma_{l,r},z=z^{l,r}(x)} |H^{l,r}(z)||\partial_x z_{l,r}(x)| \frac 1 {d(E,\CT_\delta^c,\gamma)}
\\ &\leq C(V) h \Sup{z\in\tilde\gamma_{l,r}} |v(z)|
 \frac 1 {d(E,\CT_\delta^c,\gamma)^3}
 \end{split}
 \end{gather}

\begin{gather}\begin{split}
\Sup{z\in\tilde\gamma_{l,r}}|J^{l,r}(v)(z)| &\leq  \Sup{z\in\tilde\gamma_{l,r}} |v(z)|
\Sup{z\in\tilde\gamma_{l,r},z=z^{l,r}(x)} |H^{l,r}(z) |\partial_x z_{l,r}(x)| \\&
\leq   \Sup{z\in\tilde\gamma_{l,r}} |v(z)|
\Sup{z\in\tilde\gamma_{l,r},z=z^{l,r}(x)} |H^{l,r}(z)||\partial_x z_{l,r}(x)| 
\\ &\leq C(V)  \Sup{z\in\tilde\gamma_{l,r}} |v(z)|
 \frac 1 {d(E,\CT_\delta^c,\gamma)^2}.  \end{split}\end{gather}
Similarly we obtain
\begin{gather}\begin{split}
\Sup{z\in\tilde\gamma_{l,r}}|I^{l,r}_-(v)(z)|
&\leq  h \Sup{z\in\tilde\gamma_{l,r}} |v(z)|
\Sup{z\in\tilde\gamma_{l,r},z=z^{l,r}(x)} |H^{l,r}(z)| |\partial_x z_{l,r}(x)| \frac 1 {d(E,\CT_\delta^c,\gamma)}
\\& \leq C(V) h \Sup{z\in\tilde\gamma_{l,r}} |v(z)|
 \frac 1 {d(E,\CT_\delta^c,\gamma)^3}
 \end{split}\end{gather}
 We denote by $\lfloor x\rfloor$ the integer part of $x$.
Hence on a progressive path $\gamma$ we have 
$$
|W_{n,\pm}^{l,r}|_\infty \leq \frac{ h^{\Entier\frac {n+1} 2\rfloor}} {d(E,\CT_\delta^c,\gamma)^{3n}},
$$ 
which means that
\begin{equation}\label{vingh}
\begin{array}{l}
W^l_{+}(-1)=1,\quad W^l_{+}(\beta)=1+O(\frac h {d(E,\CT_\delta^c,\gamma)^{6}}),\\
W^l_{-}(\beta)=1,\quad W^l_{-}(-1)=1+O(\frac h {d(E,\CT_\delta^c,\gamma)^{6}}),\\
W^r_{+}(\beta)=1,\quad W^r_{+}(1)=1+O(\frac h {d(E,\CT_\delta^c,\gamma)^{6}}),\\
W^r_{-}(1)=1,\quad W^r_{-}(\beta)=1+O(\frac h
{d(E,\CT_\delta^c,\gamma)^{6}}).\end{array}
\end{equation}
The function $w^l_\pm$ form a basis of  solutions of the equation $H_\delta
w=0$, and similarly $w^r_\pm$ is a basis of  solutions of the equation
$H_\delta w=0$.
A computation gives 
\begin{equation}\nonumber
\frac {dw^{l,r}_\pm}{dx}(x)=\pm {\partial_x z_{l,r}(x) }^{\frac 1 2} e^{\pm \frac {z_{l,r}(x)} h} \frac 1 h
\sum_{n=0}^\infty (-1)^n W^{l,r}_n(z_{l,r}(x)),
\end{equation}
and so
\begin{equation}\nonumber
\frac {dw^{l,r}_\pm}{dx}(x)=\pm  {\partial_x z_{l,r}(x)}^{\frac 1 2} e^{\pm
  \frac {z_{l,r}(x)} h}\frac 1 h(1+O(\frac h {d(E,\CT_\delta^c,\gamma)^{6}}
)).
\end{equation}
Now compute the determinant;  we  shall write  $\O=\O(\frac h
{d(E,\CT_\delta^c,\gamma)^{6}})$ for symplicity. 
\begin{gather}\begin{split}
&\det\left(\begin{array}{llll} 
w^l_-(-1)&w^l_+(-1)&0&0\\
 w^l_-(\beta)&w^l_+(\beta)&-w^r_-(\beta)&-w^r_+(\beta))\\
{w^l}'_-(\beta)&{w^l}'_+(\beta)&-{w^r}'_-(\beta)&-{w^r}'_+(\beta)\\
0&0&w^r_-(1)&w^r_+(1)\end{array}\right)=
\\&
\frac 1 h \partial_x z_r^{-\frac 1 2}(1)\partial_x z_r^{-\frac 1 2}(\beta)
\partial_x z_l^{-\frac 1 2}(\beta)\partial_x z_l^{-\frac 1 2}(-1)
e^{-\frac {z_l(-1)} h } e^{\frac {z_r(1)} h }\times
\\&\left|\begin{array}{llll} 
1+\O&
 e^{2\frac {z_l(-1)} h } &0&0\\
1 & 1+\O
 &-1+\O&-1 \\
- \partial_x z_l(\beta) (1+\O )   &
 \partial_x z_l(\beta)(1+\O)&
\partial_x z_r(\beta)(1+\O)&
-\partial_x z_r(\beta)(1+\O) \\
0&0& e^{-2\frac {z_r(1)} h}&1+\O\end{array}\right|
\end{split}\end{gather}
So the sign of the  determinant is given by the sign of
\begin{gather}\begin{split}
& \partial_x z_l^{-\frac 1 2}(-1)\partial_x z_r^{-\frac 1 2}(\beta)
\partial_x z_l^{-\frac 1 2}(\beta)\partial_x z_r^{-\frac 1 2}(1)\times
\\&\quad\quad
(\det\left(\begin{array}{llll} 
1& 0 &0&0\\
0 & 1 &-1&0 \\
0&\partial_x z_l(\beta)&\partial_x z_r(\beta)&0\\
0&0&0&1\end{array}\right)
+\O(\frac h {d(E,\CT_\delta^c,\gamma)^{6}}))
\\&=\partial_x z_l^{-\frac 1 2}(-1)\partial_x z_r^{-\frac 1 2}(\beta)
\partial_x z_l^{-\frac 1 2}(\beta)\partial_x z_r^{-\frac 1 2}(1)\times
\\& (\partial_x z_r(\beta)+\partial_x z_l(\beta)+\O(\frac h {d(E,\CT_\delta^c,\gamma)^{6}}) )
\\&=
{(V(-1)-i\delta-E)}^{-\frac 1 4} 
{(V(\beta)-i\delta-E)}^{-\frac 1 4}\\&\times 
{(V(\beta)+i\delta-E)}^{-\frac 1 4}
{(V(1)+i\delta-E)}^{-\frac 1 4}
\\&\times
(\sqrt{V(\beta)+i\delta-E}+\sqrt{V(\beta)-i\delta-E}+\O(\frac h {d(E,\CT_\delta^c,\gamma)^{6}}) )
\end{split}\end{gather}
( we have $\re(\sqrt{V(\beta)+i\delta-E}+\sqrt{V(\beta)-i\delta-E})\geq 2 d(E,\CT_\delta^c,\gamma)$.)

The determinant is non zero as soon as $\frac h {d(E,\CT_\delta^c,\gamma)^{6}}$ is small and
 also small 
compared to $d(E,\CT_\delta^c,\gamma)$ (and this explain the exponent $7$).\par

\cqfd

%

\section{spectrum for large $E$ and $\delta=0$}
\label{large}

In this section we prove three results: the first, Theorem \ref{reci}, gives the condition for $E_0$ to belong to $\lim_h \Sp(H^h)$.
The second,  Theorem \ref{curve}, describes  $\lim_h \Sp(H^h)$  as  a curve and  gives its asymptotics near $\infty$. 
The last, Theorem \ref{quanti}, describes the eigenvalues of $H^h$ which are  near to $\lim_h \Sp(H^h)$. 
In the following we let denote by $Y$ a primitive of the potential $V$.
In this section, we assume that $V$ verifies the following hypothesis:
\textbf{(H1)} For any $\Omega\subset \C$, if $V^{-1}(\Omega)$ bounded  then $\Omega$ is relatively compact.

\textbf{(H2)} If  $E$ is  large enough,but  with small imaginary part, then  $1$ and $-1$ belong to the same Stokes region. We recall that this means that one can find a path going from $-1$ to $ 1$ that does not intersect the Stokes lines issuing from the turning points.
\medskip
\noindent
Notice that, under the assumption \textbf{(H2)}, $E\in \CT^c_0$ is  equivalent to $\re(\int_{-1}^1\sqrt {V-E}\,dx)=0$. 
\begin{theorem}\label{reci}
For $E_0$ large  with small enough imaginary part,  then  \linebreak
$\re(\int_{-1}^1\sqrt {V-E_0}\,dx)=0$ if and only if  $E_0\in\lim_h\Sp(H^h_0)$.

Moreover there exists $C$ small enough and  $E_h$  in the spectrum of  $H^h_0$ which satisfies   $|\sqrt{ E_h}-\sqrt{ E_0}|\leq C h$.
\end{theorem}

{\sl Proof of Theorem \ref{reci}:}
The reverse implication  follows directly from Theorem  \ref{inclus}, so we prove the direct implication.

There exists two WKB solutions $w_\pm$ of $H^h_0-E=0$ of the form 
\begin{equation}
w^\pm(x) ={\partial_x z_0(x)}^{-\frac 1 2}e^{\pm z_0(x)/h}W_\pm(z_0(x)),
\end{equation}
$$W_\pm=\sum_{j=0} W_{n,\pm}\quad    W_{0,\pm}(-1)=1, \quad W_{n,\pm}(-1)=0 \mbox{ for } n>0.$$
The complex number $E$ belongs to the spectrum of $H^h_0$ if and only if 
$$
e^{\frac 2 h  z_0(1)}=\frac {W_+(1)}{W_-(1)},
$$
which is equivalent to the  existence of  $k\in\Z$ such that
$$
\int_{-1}^1 \sqrt{V(x)-E}\,dx-\frac h 2 \ln(W_+(1)) +\frac h 2 \ln(W_-(1)) -ihk\pi=0.
$$
Let us denote
 $$
 f(E)=\int_{-1}^1 \sqrt{V(x)-E}\,dx-\im\int_{-1}^1 \sqrt{V(x)-E_0}\,dx,
 $$
 and 
  $$
  k_0=\lfloor\frac {\im\int_{-1}^1 \sqrt{V(x)-E_0}\,dx} h\rfloor
  $$
 and also
$$
g(E)=\int_{-1}^1 \sqrt{V(x)-E}\,dx-\frac h 2 \ln(W_+(1)) +\frac h 2 \ln(W_-(1)) -i hk_0\pi.
$$
We have  $f(E_0)=0$. We want to apply  Rouch\'e's theorem,
but to do so we must give  upper and lower bounds for  $W_\pm(1)$ and an estimate on 
$f'(E_0)=-\frac 1 2 \int_{-1}^1\frac 1 {\sqrt {V(x)-E_0}} \,dx$.

By the hypothesis on the geometry of the Stokes lines near $-1$ and $1$ there exists a path $\gamma$  which links $-1$ to $1$  on which $\re\int_{-1}^{t}\sqrt{V(x)-E_0}\,dx=0$ for all $t\in\gamma$. Moreover $V(\gamma)$ is bounded uniformly with respect to $E_0$, and by assumption \textbf{(H1)}, $\gamma$ is also bounded uniformly with respect to $E_0$. 

We have 
$$
\frac 1 {\sqrt {V(x)-E_0}}=\frac 1 {i\sqrt E_0} (1+\O(\frac {sup_\gamma V} {E_0})),
$$
therefore
$$
|f(E)-f(E_0)|\geq \frac 1 {2\sqrt{| E_0|}} |(E-E_0)|,
$$ 
$$
\sup_{\gamma} |H(z(x))|\leq C_1 \frac {\Sup{\gamma} |V'(x)| } {|E_0|^{\frac 3 2}},
$$
and
$$
\sup_{b,b'\in\gamma} |\re\int_b^{b'} \sqrt{V(t)-E}\,\,dt|\leq C_2 \frac {E-E_0} {\sqrt {|E_0|}},
$$
 for any $E$ and $\a_0$ such that $|E-E_0|\leq \po (E_0)$  and  $|E_0|>C_0$ .  
Using the expression of
$W_{n,+}(z)$ as a Volterra integral,
we obtain the estimates
$$
|W_{n,\pm}(z)|\leq \exp\big\{\frac 2 h \sup_{b,b'\in\gamma} |\re\int_b^{b'} \sqrt{V(t)-E}\,dt| (n+2) \big\} \sup_\gamma |H|^n \frac 1 {n!}  .
$$
This gives
$$
|W_{n,\pm}(z)|\leq \exp\big\{\frac   {C_2} h \frac {|E-E_0|} {\sqrt {|E_0|}} {n+2} \big\} \frac { \hat C_1^n }{|E_0|^{\frac {3n} 2} n!}
$$
for $|E-E_0|\leq h_{0}(E_0)$ and  $|E_0|>C_0$ .
So we obtain 
$$
W_\pm(1)=1+  \frac { \hat C_1 }{|E_0|^{\frac 3 2} }
 \exp \big\{ \frac {2 C_2} h \frac {|E-E_0|} {\sqrt{ |E_0|}}\big\}
 \exp\big\{ \frac {\hat C_1}{|E_0|^{\frac 3 2} }e^{\frac  {2 C_2} h } \frac {|E-E_0|} {\sqrt{ |E_0|}     }\big\}  .
$$
This gives the estimate
$$|W_\pm(1)-1|\leq   \frac { \tilde C_1 }{|E_0|^{\frac 3 2} } e^{C_2 C_3} e^{ \hat C_1  {|E_0|^{-\frac 3 2}}
e^{  C_2 C_3 }}$$ 
for $E_0>C_0$  and $|E-E_0|\leq C_3 h\sqrt{|E_0|}$ with $h$ small enough. We also have
\begin{equation}\label{vingt}
 |\ln W_+(1)|+|\ln W_-(1)|\leq  \frac {\check C_1 } {|E_0|^{\frac 3 2}} \end{equation}
for $E_0>C(C_3)$, and $|E-E_0|\leq C_3 h\sqrt{|E_0|}$ with $h$ enough small,
and  for the functions $f$ and $g$ we get
$$
|f(E)-g(E)|\leq h +\frac h 2 (|\ln W_+(1)|+|\ln W_-(1)|).
$$
Thus if $|E-E_0|= C_3 h\sqrt{|E_0|}$, for some large constant $C_{3}>0$,  we have
\begin{equation}
\label{titus}
|f(E)-g(E)|\leq (1+C_4 |E_0|^{-\frac 3 2} )  h,\quad \quad |f(E)|\geq \frac {C_3}2h.
\end{equation}
Finally, if $E_0$ is large enough, all the assumptions of Rouch\'e's theorem are fulfilled, and  we get the result.\cqfd

\begin{theorem}\sl
\label{curve}
For any fixed $a$  large enough, the equation
$$
\re S_{-1,1}(E)_{\vert_{E=a+ib}}=\re\int_{-1}^1\sqrt{V(x)-(a+ib)}\,dx=0
$$
has a  unique solution $b(a)$. Moreover 
$b(a)=  \frac 1 {2} \im(Y(1)-Y(-1))+\O( \frac 1 {a} )$.
\end{theorem}

{\sl Proof :} First we prove that if  $b$ is such that
$\re\int_{-1}^1\sqrt{V(x)-a-ib}\,dx=0$,  and
$b=\po(a)$ then 
\begin{equation}
b= i \frac 1 2 \im(Y(1)-Y(-1))+\O( a^{-1} ).
\label{etape1}
\end{equation}
Indeed if we  denote by $E=a+ib$
$$
\varphi(E,\alpha,y)= \int_{\alpha}^y\sqrt{V(x)-E}\,dx,
$$ 
then, uniformly for $y$ in a compact set,
$$
\varphi(E,\alpha,y)=i\sqrt E (y-\alpha)-i\frac 1 {2\sqrt E} (W(y)-W(\alpha))+O(E^{-\frac 3 2}).
$$
Since $\re S_{-1,1}(E)=0$, we have
$$
\re\big(i2\sqrt E -i\frac 1 {2\sqrt E} (Y(1)-Y(-1))\big)=\O(E^{-\frac 3 2}).
$$
Writing $\sqrt E=c+i\tilde b$, then $\tilde b=\po(c)$ and
$$
-2\tilde b +\frac 1 {2c} \im(Y(1)-Y(-1))+O(\frac {\tilde b} {c^2})=\O(c^{- 3 }),
$$
or equivalently
$$
\tilde b =\frac 1 {4c} \im(Y(1)-Y(-1))+o(\frac 1 c).
$$
This gives $b= \frac 1 {2} \im(Y(1)-Y(-1))+o(1)$, and  $b= O(1)$. Therefore
$\tilde b=O(\frac 1 c)$ and recycling through this argument,
 $$\tilde b =\frac 1 {4c} \im(Y(1)-Y(-1))+O(\frac 1 {c^3}),
 $$
 which gives (\ref{etape1}).

The existence part of the theorem is straightforward:  set
$$
\varphi(c,\hat b):=\re\int_{-1}^1\sqrt{cV(x)-1-i\hat b}\,dx=c\re\int_{-1}^1\sqrt{V(x)-E}\,dx
$$
where $E=\frac 1 c (1+i\hat b)$. Since
$\partial_b\varphi(0,0)=-2$ and
$\varphi(0,0)=0$,
the implicit function theorem applies.

Finally, the uniqueness follows from the fact that the map
$\phi:b\mapsto\re\int_{-1}^1\sqrt{V(x)-a-ib}\,dx$ is injective for $a$ large enough. Indeed 
$$
|\sqrt{V(x)-a-ib}-\sqrt{V(x)-a-ib'}|\leq C \sqrt{|b-b'|},
$$
for a suitable branch of the square root.\cqfd.

\begin{remark}\sl  We can compute more terms of the asymptotic expansion than in (\ref{etape1}). Indeed,
from
\begin{gather}\begin{split}
\varphi(E,-1,1)=i2 \sqrt E -i\frac 1 {2\sqrt E} (Y(1)-Y(-1))+ i E^{-\frac 3 2} \int_{-1}^1 V^2(x)\,dx
\\-iE^{-\frac 5 2} \int_{-1}^1 V^3(x)\,dx +O(E^{-\frac 7 2}).
\end{split}\end{gather}
we get
$$b(a)=   \frac 1 {2} \im(Y(1)-Y(-1)) +\frac 3 {8a^2} \im(\int_{-1}^{1} V^3 )+\O( \frac 1 {a^3} ).$$
\end{remark}

Looking more carefully at the quantization rules, we obtain the asymptotics of the eigenvalues.

\begin{theorem}
\label{quanti}
Let $E_0$ be a  solution of  $\re(\int_{-1}^1\sqrt {V-E_0}\,dx)=0$ with
$|E_0|$ large and $\im(E_0)$ small.
If $E_h$ is in the spectrum of  $H^h_0$ and satisfies  $|\sqrt{ E_h}-\sqrt{ E_0}|\leq C h$, provided $C$ is small enough then $E_{h}$
satisfies:
$$ E= (\frac {\pi hk} 2)^2 +
 \frac {(Y(1)-Y(-1))} {2} + \frac {(Y(1)-Y(-1))^2} {{(2hk\pi)}^2}   +  \O( \frac 1 {(hk)^3} ).
 $$
for some $k\in \nit$.
\end{theorem}
\begin{remark}
For a real potential $V$ bounded on $[-1,1]$ and for any $E_0\in\rit $ large we get  $\re(\int_{-1}^1\sqrt {V-E_0}\,dx)=0$.
So Theorem \ref{quanti} can be read as:
If $E_h$ is large and satisfies $E_h\in\Sp(H^h_0)$ then $E_{h}$
satisfies:
$$ E= (\frac {\pi hk} 2)^2 +
 \frac {(Y(1)-Y(-1))} {2} + \frac {(Y(1)-Y(-1))^2} {{(2hk\pi)}^2}   +  \O( \frac 1 {(hk)^3} ).
 $$
for some $k\in \nit$.
\end{remark}
{\sl Proof :}
$$\int_{-1}^1\sqrt {V-E}\,dx=i2\sqrt E-i\frac 1 {2\sqrt E} (Y(1)-Y(-1))+\O(E^{-\frac 3 2}).$$
Then $E\in\Sp(H^h_\delta)$ if and only if there exists $k\in\Z$ such that
$$\int_{-1}^1 \sqrt{V(x)-E}\,dx-\frac h 2 \ln(W_+(1)) +\frac h 2 \ln(W_-(1)) -ihk\pi=0$$
For $E_0$ big enough,  $|E-E_0|\leq C_3 h\sqrt{E_0}$, $C_3$ is as in
(\ref{titus}) and  $h$  small, we get  
$$ |\ln W_+(1)|+|\ln W_-(1)|\leq  \check{ C} |E_0|^{-\frac 3 2} $$
So $E\in\Sp(H^h_\delta)$ if and only if there exist $k\in\Z$ such that
$$i2\sqrt E-i\frac 1 {2\sqrt E} (Y(1)-Y(-1)) -ihk\pi=\O(E^{-\frac 3 2})$$
In particular taking $\sqrt E=c+i\tilde b$  the real and imaginary parts give 
$$2 c - hk\pi=\O(c^{-1})\quad c =\frac 1 2 hk\pi+\O((kh)^{-1})$$
$$2\tilde b=\frac 1 {2c} \im(Y(1)-Y(-1))+\O( c^{-3} )$$
Using this again gives
$$c =\frac 1 2 hk\pi+\frac 1 {2hk\pi} \re(Y(1)-Y(-1))+\O((kh)^{-3})$$
\cqfd
\medskip
\begin{remark}
One can also treat the case $P^h$ on all of $\rit$.  Then a condition for $E$ to belong to 
 $\lim_h \Sp(H^h)$ is that  there exists a progressive path joining two  turning points. 
 \end{remark}
%

\section{Large spectrum of  a perturbation of H}
\label{largeper}

In this section, we prove two theorems: Theorem \ref{reciper}  gives the condition for $E_0$ to belong to $\lim_h \Sp(H^h_\delta)$ and Theorem  \ref{curveper} gives a description  of $\lim_h \Sp(H^h_\delta )$ as
the union of two curves, the asymptotics  near $\infty$ of which are made explicit. 

We assume that $V$ verifies the  hypothesis  \textbf{(H1)} and in addition~:  
\par
\textbf{(H3)} If  $E$ is  large enough, with small imaginary part then  $1$,  $-1$  and $\beta$ belong to the same Stokes region.\par
\medskip
We assume $\beta\neq \pm 1$.
Recall that \textbf{(H3)}  means that one can find  paths going from $ -1$ to
$ \beta$  and from $ 1$ to $\beta$  which do not intersect the Stokes lines issuing from the turning points.
\medskip
 With the hypothesis \textbf{(H3)}, the condition $E\in \CT^c_\delta$ is then equivalent to \par\noindent
either $\re(\int_{-1}^\beta\sqrt {V-i\delta-E}\,dx)=0$  
or $\re(\int_{\beta}^{1}\sqrt {V+i\delta-E}\,dx)=0$.  

\begin{theorem}\label{reciper}
For $E_0$ large enough,  with small enough imaginary part,
 then $E_0$ satisfies $\re(\int_{-1}^\beta\sqrt {V-E_0}\,dx)=0$ or $\re(\int_{\beta}^1\sqrt {V-E_0}\,dx)=0$ if and only if $E_0\in\lim_h\Sp(H^h_\delta)$
 \end{theorem}

\begin{theorem}\label{curveper}
For any fixed $a$ large enough, the equation 
$$\re(S^\delta_{\pm1,\beta}(E))|_{E=a+ib}:=\re\int_{\pm1}^\beta \sqrt{V(x)\pm i\delta-a-ib}\,dx=0$$
 has a unique solution $b(a)$, and this solution  satisfies
$$ b(a)= i \frac 1 {\beta-\pm 1} \im(Y(\beta)-Y(\pm 1)) \pm i \delta +O( \frac 1 a ).$$
\end{theorem}

{\sl Proof :}
The proof is the same as in Theorem \ref{curve} where the potential is $V\pm i\delta$  in each side and we just need to change the interval for the integral. \cqfd

{\sl Proof of Theorem \ref{reciper}:}
The reverse implication  is already proved in Theorem \ref{inclus}. Let us
prove the direct implication as in the proof of Theorem \ref{reci}.
There exist four WKB solutions $w^{l,r}_\pm$ of $H^h_\delta-E=0$  on the interval $[-1,\beta]$ or $[\beta,1]$
\begin{equation}\begin{array}{l}
w^l_\pm =({\partial_xz_l}(x))^{-\frac 1 2}e^{ \pm z_l(x)/h} W^l_\pm (z_l(x)),\\
w^r_\pm =({\partial_x z_r}(x))^{-\frac 1 2}e^{ \pm z_r(x)/h} W^r_\pm (z_r(x)),\end{array}
\label{wkbsolld}
\end{equation}
with initial data
\begin{equation}
\left\{
\begin{array}{ll}
W^{l,r}_{0,\pm}=1,&\\[8pt]
W^{l}_{n,+}(z_l(\beta))=0,\ W^{r}_{n,-}(z_r(\beta))=0, & n\geq 1,\\[8pt]
W^{l}_{n,-}(z_l(\beta))=0,\ W^{r}_{n,+}(z_r(\beta))=0, & n\geq 1.
\end{array}
\right .
\end{equation}

Recall that $E\in\Sp(H^h_\delta)$ if and only if 
\begin{gather}\begin{split}\label{detcal}
&0=\det\left(\begin{array}{llll} 
w^l_-(-1)&w^l_+(-1)&0&0\\
 w^l_-(\beta)&w^l_+(\beta)&-w^r_-(\beta)&-w^r_+(\beta)\\
{w^l}'_-(\beta)&{w^l}'_+(\beta)&-{w^r}'_-(\beta)&-{w^r}'_+(\beta)\\
0&0&w^r_-(1)&w^r_+(1)\end{array}\right)=
\\&
\frac 1 h \partial_x z_r^{-\frac 1 2}(1)\partial_x z_l^{-\frac 1 2}(-1)\times
\\& \det\left(\begin{array}{llll} 
{\scriptscriptstyle e^{-\frac{ z_l(-1)} h} W^l_- (z_l(-1))}
 &{\scriptscriptstyle e^{\frac {z_l(-1)} h }W^l_+ (z_l(-1))}&0&0\\
1 & 1 &-\frac { \partial_x z_l^{\frac 1 2}(\beta)}{\partial_x z_r^{\frac 1 2}(\beta)} 
&-  \frac {\partial_x z_l^{\frac 1 2}(\beta)}{\partial_x z_r^{\frac 1 2}(\beta)} \\
-1 &1& \frac {\partial_x z_r^{\frac 1 2}(\beta)} {\partial_x z_l^{\frac 1 2}(\beta)}
&-\frac {\partial_x z_r^{\frac 1 2}(\beta)}{ \partial_x z_l^{\frac 1 2}(\beta)}\\
0&0&{\scriptscriptstyle e^{-\frac {z_r(1)} h}W^r_- (z_r(1))}&
{\scriptscriptstyle e^{\frac {z_r(1)} h}  W^r_+ (z_r(1))} \end{array}\right)
\end{split}\end{gather}
To make the computation,  simplify the notation by setting 
$$t=\frac { \partial_x z_l^{\frac 1 2}(\beta)}{\partial_x z_r^{\frac 1
    2}(\beta)} ,\;y=e^{\frac {z_r(1)} h}\; x= e^{\frac {z_l(-1)} h}$$
then (\ref{detcal}) becomes
 \begin{gather}\begin{split}
&\det\left(\begin{array}{llll} 
x^{-1} W^l_- (z_l(-1)) & xW^l_+ (z_l(-1))&0&0\\
1 & 1 &-t &- t \\
-1 &1&  t^{-1}  &- t^{-1} \\
0&0& y W^r_- (z_r(1))& y^{-1}  W^r_+ (z_r(1)) \end{array}\right)=0,
\end{split}\end{gather}
which gives
$$ x^2 \frac {W^l_+ (z_l(-1))}{W^l_- (z_l(-1))}= \big (1+y^2  \frac {W^r_-
  (z_r(1))}{W^r_+ (z_r(1))} \frac {1-t^2}{1+t^2} \big) \big({\frac {1-t^2}{1+t^2} +y^2  \frac {W^r_- (z_r(1))}{W^r_+ (z_r(1))} }^{-1}\big).$$
We remark that 
$$|\big(1-\frac {\partial_x z_l(\beta)}{\partial_x z_r(\beta)} \big) \big({1+\frac {
    \partial_x z_l(\beta)}{\partial_x z_r(\beta)}}^{-1}\big)|\leq C_8 \sqrt \delta.$$ 
So the condition for $E$ to be in the spectrum of $H_\delta^h$ is 
\begin{gather}\begin{split}
& \big(e^{2\frac {z_l(-1)} h}  \frac {W^l_+ (z_l(-1))}{W^l_- (z_l(-1))}-
\frac {(\partial_x z_r(\beta)-\partial_x z_l(\beta))}{(\partial_x
  z_r(\beta)+\partial_x z_l(\beta))} \big)\times
\\&
\big(e^{-2\frac {z_r(1)} h}  \frac {W^r_- (z_r(1))}{W^r_+ (z_r(1))}+\frac
{(\partial_x z_r(\beta)-\partial_x z_l(\beta))}{(\partial_x
  z_r(\beta)+\partial_x z_l(\beta))}\big)
\\&=\big(1-\frac {(\partial_x z_r(\beta)-\partial_x z_l(\beta))^2}{(\partial_x
  z_r(\beta)+\partial_x z_l(\beta))^2}\big).\end{split}\end{gather}

We write this condition by taking logs as 
\begin{gather}\begin{split}\label{titi}
&\ln\big( e^{2\frac {z_l(-1)} h}  \frac {W^l_+ (z_l(-1))}{W^l_- (z_l(-1))}-
\frac {(\partial_x z_r(\beta)-\partial_x z_l(\beta))}{(\partial_x
  z_r(\beta)+\partial_x z_l(\beta))}\big)
\\&+\ln\big(e^{-2\frac {z_r(1)} h}  \frac {W^r_- (z_r(1))}{W^r_+
  (z_r(1))}+\frac {(\partial_x z_r(\beta)-\partial_x z_l(\beta))}{(\partial_x
  z_r(\beta)+\partial_x z_l(\beta))}\big)\\
&-\ln\big (1-\frac {(\partial_x z_r(\beta)-\partial_x z_l(\beta))^2}{(\partial_x
  z_r(\beta)+\partial_x z_l(\beta))^2}\big)=2ik\pi,\quad k\in\zit.
\end{split}\end{gather}
\underline{Step1} 
First assume that 
$$\re(\int_{-1}^\beta\sqrt {V-i\delta-E_0}\,dx)=0\;\mbox{ and }\;  \re(\int_{\beta}^1\sqrt {V+i\delta-E_0}\,dx)\neq0.$$
 (The case $\re(\int_{-1}^\beta\sqrt {V-i\delta-E_0}\,dx)\neq 0$ and $\re(\int_{\beta}^1\sqrt {V+i\delta-E_0}\,dx)=0$ is treated the same way. )\par
As  (\ref{vingt}), we get  
\begin{equation}\label{auauu} |\ln W^l_+(z_l(-1))|+|\ln W^l_-(z_l(-1))| \leq   C_1 |E_0|^{-\frac 3 2} 
\end{equation}
 for $|E-E_0|\leq C_3 h\sqrt{|E_0|}$ and $E_0>C$ big, $h$ small. 
As in (\ref{vingh}), we get
\begin{equation} \label{auau} |\ln W^r_+ (z_r(1))|+|\ln W^r_- (z_r(1))|\leq  C h \end{equation}
with $C$ uniform in $E_0$,  for $|E-E_0|\leq C_3 h\sqrt{|E_0|}$ 
and 
\begin{equation}\label{auu}| e^{-2\frac {z_r(1)} h}|\leq e^{-\frac {C} h}.\end{equation}
Rewrite the previous condition as 
\begin{gather}\begin{split}
&g_{1,k}(E)= 2 z_l(-1) +h \ln( \frac {W^l_+ (z_l(-1))}{W^l_- (z_l(-1))})
\\&
+h \ln(1-e^{-2\frac {z_l(-1)} h}  \frac{W^l_- (z_l(-1))} {W^l_+
  (z_l(-1))} \frac {(\partial_x z_r(\beta)-\partial_x z_l(\beta))}{(\partial_x
  z_r(\beta)+\partial_x z_l(\beta))})
\\&+h \ln(\frac {(\partial_x z_r(\beta)-\partial_x z_l(\beta))}{(\partial_x
  z_r(\beta)+\partial_x z_l(\beta))})
\\&+ h \ln(1+
e^{-2\frac {z_r(1)} h}  \frac {W^r_- (z_r(1))}{W^r_+ (z_r(1))}\frac
{(\partial_x z_r(\beta)+\partial_x z_l(\beta))} {(\partial_x z_r(\beta)-\partial_x z_l(\beta))} )
\\&-h \ln (1-\frac {(\partial_x z_r(\beta)-\partial_x
  z_l(\beta))^2}{(\partial_x z_r(\beta)+\partial_x z_l(\beta))^2})-2ikh\pi=0.
\end{split}\end{gather}
Setting 
$$f_1(E)=2 z_l(-1) +2\Im(\int_{-1}^{\beta} \sqrt{ V(x)-i\delta-E_0} \,dx)$$ and 
$$k_0=\Entier  \Im(\int_{-1}^{\beta} \sqrt{ V(x)-i\delta-E_0} \,dx)h^{-1}\rfloor,$$ we have $f_1(E_0)=0$.
Using  (\ref{auauu}), (\ref{auau}) and (\ref{auu}) gives
$$|g_{1,k_0}(E)-f_1(E)|\leq h(1+C\sqrt\delta+C h\sqrt\delta^{-1} +  C_1 |E_0|^{-\frac 3 2}).$$
Now recall that 
$$|f_1(E)-f_1(E_0)|\geq \frac 1 {\sqrt{| E_0|}} |(E-E_0)|.$$
We can  apply the Rouch\'e's Theorem to prove the existence of an eigenvalue
of $H^h_\delta$  for each $h$ at a distance
 $C_3(\delta) \sqrt{| E_0|} h$ of $E_0$ for $E_0$ large. 
 \medskip
\underline{Step2} 
Assume now  that 
$$\re(\int_{-1}^\beta\sqrt {V-i\delta-E_0}\,dx)=\re(\int_{\beta}^1\sqrt {V+i\delta-E_0}\,dx) =0.$$
As in  (\ref{vingt}),  
\begin{gather}\begin{split}
|\ln W^l_+ (z_l(-1))|+|\ln W^l_- (z_l(-1))|+ |\ln W^r_+ (z_r(1))|+\\
|\ln W^r_- (z_r(1))|
\leq  C_9   |E_0|^{-\frac 3 2} .\end{split}
\end{gather}
For $E_0>C$, $C$ sufficiently large, and $|E-E_0|\leq C_3 h\sqrt{|E_0|}$ with $h$ small. 
We write (\ref{titi}) as 
\begin{gather}\begin{split}
&g_{2,k}(E)=2 z_l(-1)-2 z_r(1)  +h \ln( \frac {W^l_+ (z_l(-1))}{W^l_- (z_l(-1))})\\&
+h \ln(1-e^{-2\frac {z_l(-1)} h}  \frac{W^l_- (z_l(-1))} {W^l_+ (z_l(-1))} \frac {(\partial_x z_r(\beta)-\partial_x z_l(\beta))}{(\partial_x z_r(\beta)+\partial_x z_l(\beta))})
\\&
h \ln( \frac {W^r_+ (z_r(1))}{W^r_- (z_r(1))})
+h \ln(1-e^{2\frac {z_r(1)} h}  \frac{W^r_- (z_r(1))} {W^r_+ (z_r(1))}
\frac {(\partial_x z_r(\beta)-\partial_x z_l(\beta))}{(\partial_x
  z_r(\beta)+\partial_x z_l(\beta))})
\\& -h \ln (1-\frac {(\partial_x z_r(\beta)-\partial_x
  z_l(\beta))^2}{(\partial_x z_r(\beta)+\partial_x z_l(\beta))^2}-2ikh\pi=0,
\end{split}\end{gather}
and write
\begin{gather}\begin{split}
f_2(E)=2 z_l(-1)-2 z_r(1) -2\Im(\int_{\beta}^{-1} \sqrt{ V(x)-i\delta-E_0} \,dx) \\+2\Im(\int_{\beta}^{1} \sqrt{ V(x)+i\delta-E_0} \,dx).\end{split}
\end{gather}
Let 
$$k_0=\Entier \Im(\int_{\beta}^{-1} \sqrt{ V(x)-i\delta-E_0} \,dx-\int_{\beta}^{1} \sqrt{ V(x)+i\delta-E_0} \,dx)h^{-1}\rfloor.$$
Then
$$|f_2(E)-f_2(E_0)|\geq {\sqrt{| E_0|}}^{-1} |(E-E_0)|,$$
and for $|E-E_0|\leq C_3 h\sqrt{|E_0|}$
$$|f_2(E)-g_{2,k_0}(E)|\leq h( C_9 |E_0|^{-\frac 3 2} +C_8\sqrt\delta +1).$$
On the set where $|E-E_0|= C_3 h\sqrt{|E_0|}$, we have for $C_3=C_3(\delta)$ and  $E_0$ large
$$|(f_2-g_{2,k_0})(E)|\leq |f_2(E)| .$$
All the hypotheses of Rouch\'e's theorem are now satisfied and we conclude as before.\cqfd
%

\section{Application to the harmonic oscillator}
\label{harmonic}

In this section we are interested in  the Dirichlet  eigenvalues  of the operator $P^h$ on $L^2([-1,1])$,
\begin{equation}
 P^h=-h^2 \frac d {\,dx^2} +ix^2,
\label{Px^2}
\quad
\left\{
\begin{array}{ll}
P^hu=Eu,\\
u(-1)=u(1)=0.
\end{array}
\right .
\end{equation}

\begin{figure}[ht]
\begin{center}
\epsfysize=6cm
\epsfbox{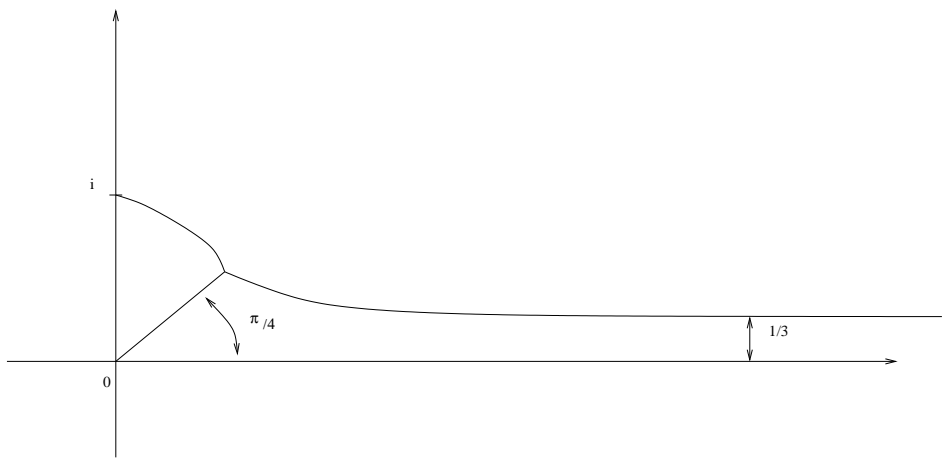}
\end{center}
\caption{Limit spectrum of $P^h$ }
\label{spectre}
\end{figure}

\begin{figure}[ht]
\begin{center}
\epsfysize=6cm
\epsfbox{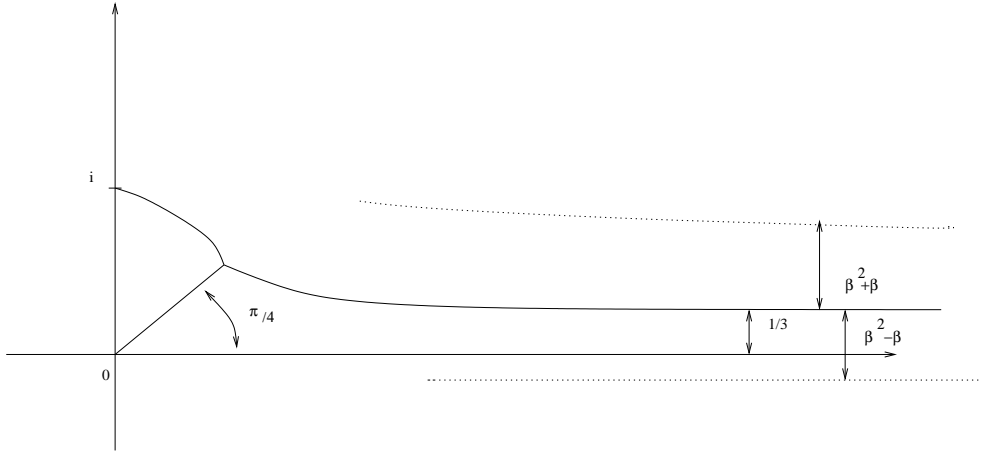}
\end{center}
\caption{Spectrum of $P^h$ in plain  lines and spectrum of $P^h_\delta$  for $\delta$ small near to $\infty$ as dotted lines.}
\label{spectreper}
\end{figure}

We obtain here only two results. Theorem \ref{forme}  describes  the shape of
the set $\lim_h\Sp(P^h)$, cf.  figure \ref{spectre}.
Theorem \ref{infiniper}  shows how the spectrum changes near $\infty$
when $P^h$ is change by a specific perturbation of size $\delta$, as illustrated in figure
 \ref{spectreper} . 

For any $E\in \C^*$ there are two turning points, $\alpha_\pm (E)=\pm(-iE)^{1/2}$, which are both simple i.e. of order 1. 

\begin{figure}[ht]
\begin{center}
\epsfysize=5cm
\epsfbox{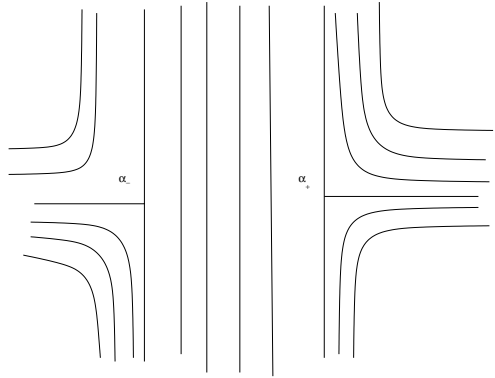}
\end{center}
\caption{Non degenerate Stokes lines}
\label{gen1}
\end{figure}

Denote by $S_{x,y}(E)=\int_x^y (it^2-E)^{1/2} \,dt$.
The function $S_{0,x}$ is even, so to simplify the computation one can use that 
$S_{-x,x}=2S_{0,x}$.
Writing
$$z(x,y,E)=\int_x^y (it^2-E)^{1/2} \,dt, $$
$$z^{har}(x,y,E)=\int_x^y (u^2-E)^{1/2} \,du,$$
and changing  coordinate $t=e^{-\frac{i\pi} 8}u$, we get 
 $$z(x,y,E)=\pm z^{har}(x e^{\frac{i\pi} 8} ,y e^{\frac{i\pi} 8} ,Ee^{-\frac{i\pi} 4}).$$
 So the Stokes lines of $-h^2\Delta+ix^2$ can be deduce from those  of the harmonic oscillator 
 $-h^2\Delta+x^2$ by a rotation by  $-\frac {\pi} 8$ about the origin,  \cite{fe}.
Using the geometry of the Stokes lines we find for some $E$ a progressive path from $-1$ to $1$. Combine with previous theorems, we obtain a partial description   of  $\CT^c$~:  
\begin{theorem} We have\par\noindent
 \begin{itemize}
 \item Suppose $E\in \C^*$ is such that 
 $
 \re S_{ \alpha_-,\alpha_+}(E)\neq 0,$  $\re S_{-1,1}(E)\neq 0,$   $\re S_{\alpha_+,1}(E)\neq 0, $  $\re S_{\alpha_-,1}(E)\neq 0,$
 then $E\in\CT$
\item  If $E\in \CT$, then for any $h$ small enough, $E$ is not an eigenvalue for $P^h$.
\item  If $E\in\lim_h\Sp(P^h)$ then $E\in\CT^c$
\end{itemize}
\label{main}
\end{theorem}
{\sl Proof:}
Using Theorem \ref{mainbis} and Theorem \ref{inclus}, it is enough to find a path  $\gamma$ tranversal to the Stokes lines from $-1$ to $1$.
We have $\re S_{\alpha_-,\alpha_+}(E)\neq 0$,  so the Stokes line issuing from $\a_+$ do not intersect the Stokes lines issuing form $\a_-$. Therefore the complex plane is divided into exactly five region delimited by the Stokes lines issuing from the turning points. 
In  Figure \ref{gen1},  we have drawn the configuration of the Stokes lines up to an analytic diffeomorphism.\par
We have $\re S_{1,\alpha_\pm}(E)\neq 0$, so that in particular $1$, (and by
symmetry $-1$) are not on any boundary of the Stokes regions. \par
 \begin{figure}[ht]
\begin{center}
\epsfysize=5cm
\epsfbox{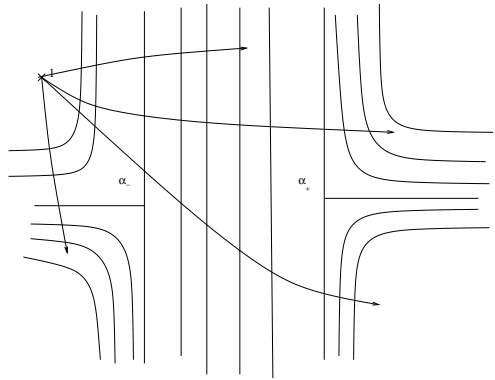}
\end{center}
\caption{Progressive path issuing from $1$}
\label{che}
\end{figure}

If $1$ and $-1$ belong to different regions, then there always exists a path $\gamma$ tranversal to the Stokes lines (see the figure \ref{che});
if on the other hand  $1$ and $-1$ belong to the same region, but we also have
$\re S_{-1,1}(E)\neq 0$, then they do not belong to the same Stokes lines and  there still exists  such a path $\gamma$.  \cqfd
\par
\begin{figure}[ht]
\begin{center}
\epsfysize=5cm
\epsfbox{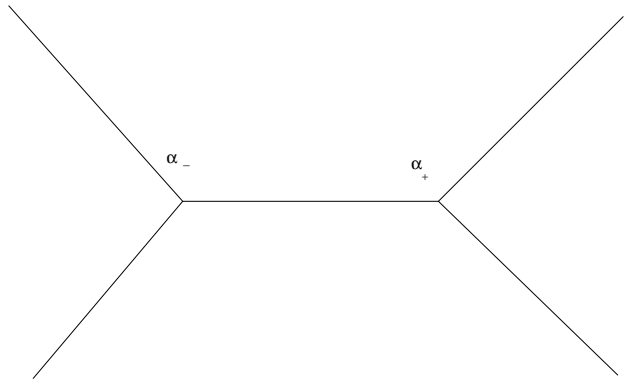}
\end{center}
\caption{Degenerate Stokes lines}
\label{gen2}
\end{figure}

\medskip
If $\re S_{\alpha_-,\alpha_+}(E)= 0$,  the Stokes configuration is
as in Figure \ref{gen2}.\par
We now picture the evolution of the Stokes lines as $E$ moves in the complex plane.\par
\begin{figure}[ht]
\begin{center}
\epsfysize=8cm
\epsfbox{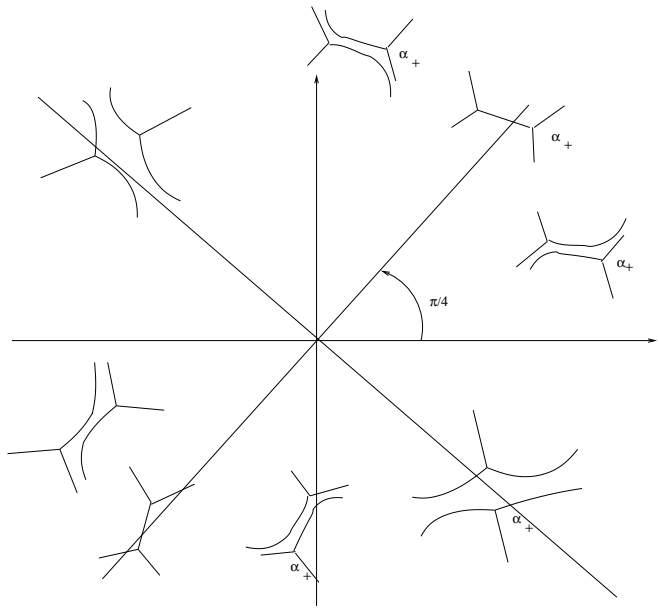}
\end{center}
\caption{Stokes line for different $E$}
\label{fullstoke}
\end{figure}
\begin{figure}[ht]
\begin{center}
\epsfysize=4cm
\epsfbox{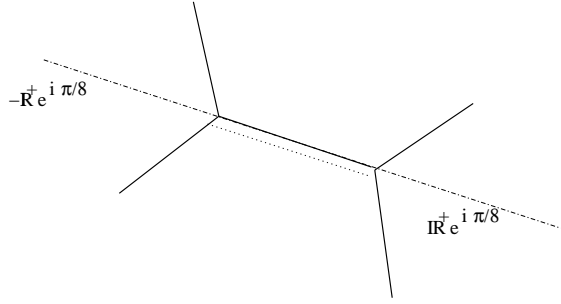}
\end{center}
\caption{branch for $E=e^{i\frac \pi 4}$, the cut is the dashed 
  line and also a Stokes line}
\label{det2}
\end{figure}
\begin{figure}[ht]
\begin{center}
\epsfysize=4cm
\epsfbox{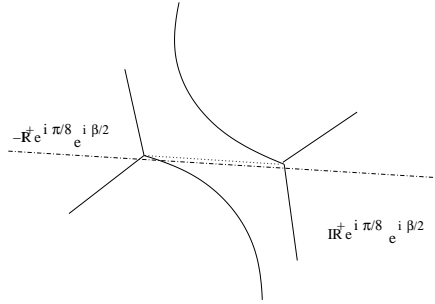}
\end{center}
\caption{branch  for $E=e^{i\frac \pi 4}e^{i\beta}$,$\beta>0$ small. The
  cut is the dashed line }
\label{secondet}
\end{figure}

In fact  Theorem \ref{main} can  be improved.  Choose  the segment $[\a_-,\a_+]$ as a cut for  $x\mapsto (V(x)-E)^{1/2}$. 
The pictures \ref{det2} and \ref{secondet} show the branch of  the square root  for different values of $E$. \par 
\par

\begin{figure}[ht]
\begin{center}
\epsfysize=3cm
\epsfbox{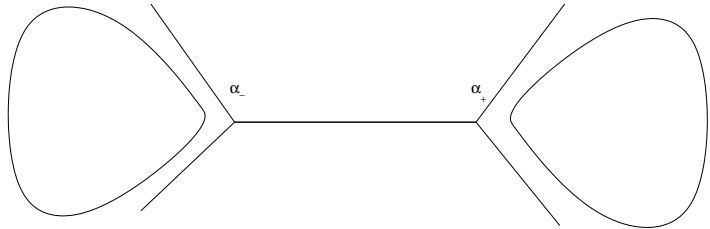}
\end{center}
\caption{The position of $1$ is in one of hatched regions}
\label{sit1}
\end{figure}

 Denote by 
 $$\Gamma_{-1,1}=\{E\in\cit;\re S_{-1,1}(E)= 0,\quad \re S_{\alpha_-,\alpha_+}(E)\neq 0\}$$
The next theorem is proved by listing all the situations where there do not exist  progressive paths from $-1$ to $1$ using the geometry of the Stokes lines.  
\pagebreak
\begin{theorem}\sl\label{main2}We have\par\noindent
 \begin{itemize}
 \item
 Suppose $E\in \C^*$  belongs to 
$$\Gamma_{\alpha_-,\alpha_+}=\{E\in\cit;  \re S_{\alpha_-,\alpha_+}(E)= 0, \re S_{\alpha_+,1}(E)\leq 0\},$$ then $E\in\CT$.
\item
Suppose $E\in \C^*$  belong to 
 $$\Gamma_{\alpha_+,1}=\{E\in\cit; \re S_{\alpha_+,1}(E)= 0; \im S_{\alpha_+,1}(E)\leq 0\},$$ then $E\in\CT$.
\item
Suppose $E\in \C^*$  belong to 
$$\Gamma_{\alpha_-,1}=\{E\in\cit; \re S_{\alpha_-,1}(E)= 0; \im S_{\alpha_-,1}(E)\geq 0\},$$ then $E\in\CT$.
 \item $E\notin \Gamma_{\alpha_+,1}\cup\Gamma_{\alpha_-,1}\cup\Gamma_{\alpha_-,\alpha_+}\cup\Gamma_{-1,1}$ then there is no progressive path form $-1$ to $1$. 
  \end{itemize}
\end{theorem}
\begin{remark}\label{stupid}
If a progressive path enters a  Stokes region crossing a Stokes line issuing from a turning point $\alpha$, then it cannot leave this region by crossing any others    Stokes lines issuing from the same turning point $\alpha$.  
\end{remark}
The proof consists in listing the cases where progressive paths do not exist.
Recall first that if  $1$ and $-1$ belong to the same region then there exists a progressive path if and only if 
$\re S_{-1,1}(E)\neq 0$. 
We get   $\Gamma_{-1,1}=\{E\in\cit; 1\mbox{ and }-1 \mbox{ belong to the same region }\re S_{-1,1}(E)= 0\}$.
Now  assuming  $1$ and $-1$ belong to different regions, we get two different figures for the Stokes lines either 
  $\re S_{\alpha_-,\alpha_+}(E)= 0$ (Figure\ref{gen2}) or  $\re S_{\alpha_-,\alpha_+}(E)\neq 0$ (Figure \ref{gen1}).
In the first case, using Remark \ref{stupid}, we see that we cannot find a path if and only if  $1$ belongs to the hatched region ( and $-1$ by symmetry to  the opposite one)(Figure \ref{sit1}).\par
\begin{figure}[ht]
\begin{center}
\epsfysize=3cm
\epsfbox{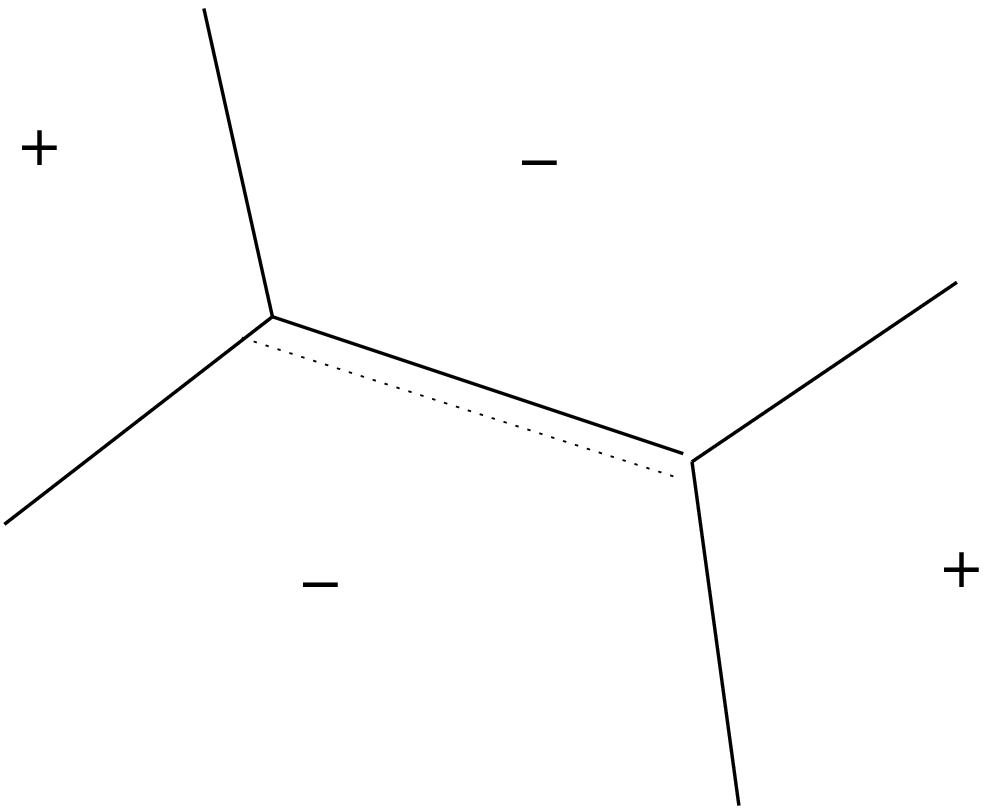}
\end{center}
\caption{ Sign of the quantity $x\to\re S_{\alpha_+,x}(e^{i\frac \pi 4})$}
\label{deter}
\end{figure}

The condition $\re S_{\alpha_-,\alpha_+}(E)= 0,\quad \re S_{\alpha_+,1}(E)>0$ corresponds to this situation, i.e.
 Figure \ref{deter} shows  the sign of the quantity $\re S_{\alpha_+,x}(e^{i\frac \pi 4}\lambda)$ depends on the position of $x$ for  $\lambda\in\rit$. \par
\begin{figure}[ht]
\begin{center}
\epsfysize=5cm
\epsfbox{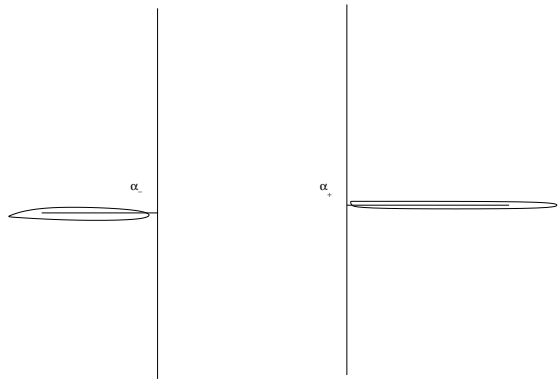}
\end{center}
\caption{The position of $1$ is in one of the dotted line, that are Stokes lines}
\label{sit2}
\end{figure}

Now  assuming  $1$ and $-1$ belong to different regions and   $\re S_{\alpha_-,\alpha_+}(E)\neq 0$ (Figure \ref{secondet}).
Using Remark \ref{stupid}, we see that we cannot find a path if and only if $1$ belong the dotted curves ( and $-1$ by symetry to the opposite one )(Figure \ref{sit2}).
The condition  $\re S_{\alpha_+,1}(E)= 0,\quad \im S_{\alpha_+,1}(E)>0$ correspond to Figure \ref{sit2}.\par
\begin{figure}[ht]
\begin{center}
\epsfysize=8cm
\epsfbox{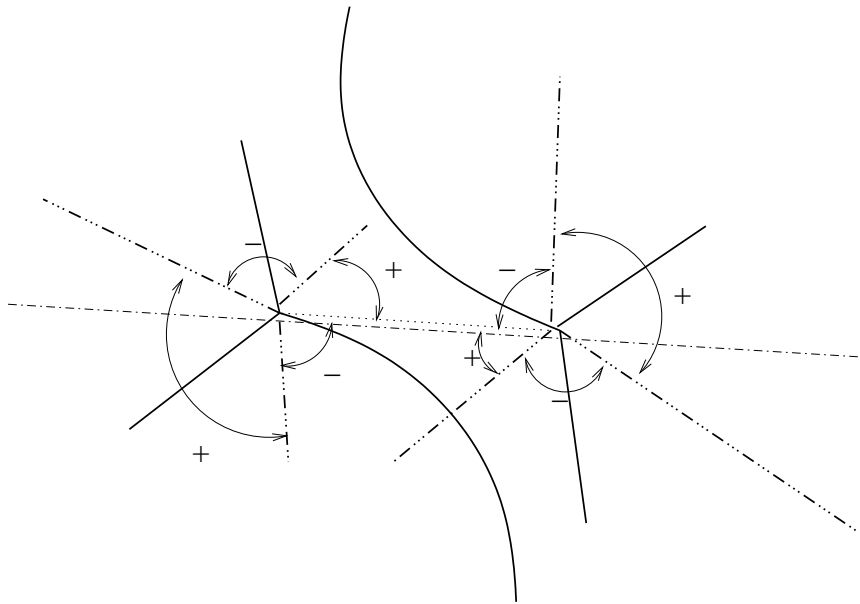}
\end{center}
\caption{ Sign of the quantity 
$\im S_{\alpha_\pm,x}(e^{i\frac \pi 4}\lambda e^{i\beta})$}
\label{seco}
\end{figure}

The condition  $\re S_{\alpha_+,1}(E)= 0,\quad \im S_{\alpha_+,1}(E)>0$ corresponds to this case.
i.e. Figure \ref{seco} shows the sign of the quantity 
$\im S_{\alpha_\pm,x}(e^{i\frac \pi 4}e^{i\beta}\lambda)$ with $\beta>0$ and small depending on $x$ near to $\alpha_\pm$.\par
Now we want to describe the set $\CT^c$.
\pagebreak
 \begin{theorem}\sl\label{forme}We have\par\noindent
 \begin{enumerate}
 \item The set $\CT^c$ is the union of three curves 
$\Gamma_{\alpha_+,1}$, $\Gamma_{\alpha_-,\alpha_+}$, $\Gamma_{-1,1}$ with only one of infinite length $\Gamma_{-1,1}$.
\item The three curves meet at a common point $\lambda_0e^{\frac {i\pi} 4}$ with $\lambda_0\in\rit^+$. 
\item The curve $\Gamma_{-1,1}$ goes from $\lambda_0e^{\frac {i\pi} 4}$ to $\rit+\frac i 3$.
\item $\Gamma_{\alpha_-,\alpha_+}=\{\lambda e^{\frac {i\pi} 4}; \lambda\in\rit, 0\leq \lambda<\lambda_0\} $.
\item The curve $\Gamma_{\alpha_+,1}$ goes from  $i$ to  $\lambda_0e^{\frac {i\pi} 4}$. 
\end{enumerate}
\end{theorem}

{\sl Proof:}
(1)~:\par
Let $E=\lambda^2 e^{i\beta}$ with $\beta\in\rit$ and $\lambda\in\rit$. We fix $\beta$ and increase $\lambda$.
We get $z(x,y,\lambda^2 e^{i\beta})=z(\frac x \lambda,\frac y \lambda,
e^{i\beta})$. The shape of the Stokes lines remains invariant up to
dilation by $\lambda$. 
Let $\beta$ be such that $\re S_{\alpha_-.\alpha_+}(e^{i\beta})\neq 0$ then under symmetry we get
 $\re S_{0,\alpha_+}(e^{i\beta})\neq 0$ , and $\re S_{\alpha_-.0}(e^{i\beta})\neq 0$. 
So if $\beta$ is such that there is no Stokes line issuing from the turning
points and going through  $0$, then there exists a neighborhood  of $0$ with no point of any  Stokes lines issuing from the turning points. 
For $\lambda$ big enough, we  get that $-1$ and $1$  belong to this neighborhood of $0$. 
Then $E$ will belong to $\CT^c$ if and only if $\re S_{-1,1}(E)\neq 0$.

\begin{figure}[ht]
\begin{center}
\epsfysize=7cm
\epsfbox{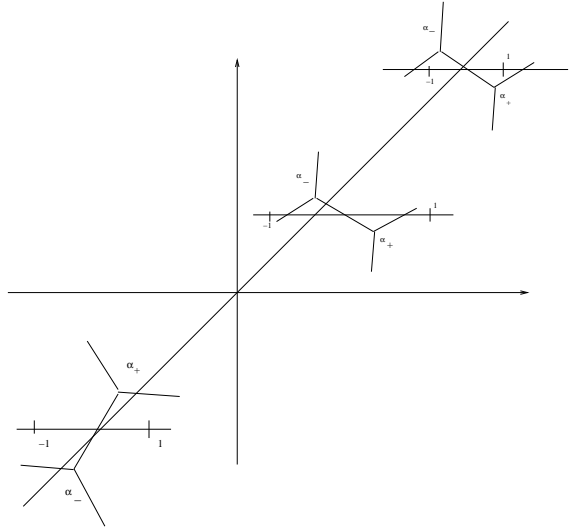}
\end{center}
\caption{Stokes line for $E=\lambda e^{\frac {i\pi} 4}$}
\label{stokesaxe}
\end{figure}

If $\beta$ is such that $\re S_{\alpha_-,\alpha_+}(e^{i\beta})= 0$   then
$\beta=\pm \frac \pi 4$ and Figure \ref{stokesaxe} and the fact that the
spectrum is included in the set $\{E\in\cit; E=\rit^++ix^2, x\in[-1,1]\}$
(i.e. the pseudo spectrum or the values of the principal symbol) shows that
$E\in\CT^c$, $E=\lambda e^{i\frac \pi 4}$, $\lambda\in\rit$,   implies $E$ bounded. 
So we have that if $E$ is big enough and belongs to $\CT^c$, then  $E\in\Gamma_{-1,1}$.

We  choose  the determination of $\alpha_\pm$ so that the set 
$$\Gamma_{\alpha_-,1}=\{E\in\cit;\, \re S_{\alpha_-,1}(E)= 0;\, \im
 S_{\alpha_-,1}(E)\geq 0\}$$
 is empty see Figure \ref{fullstoke}. 
\par
\begin{figure}[ht]
\begin{center}
\epsfysize=4cm
\epsfbox{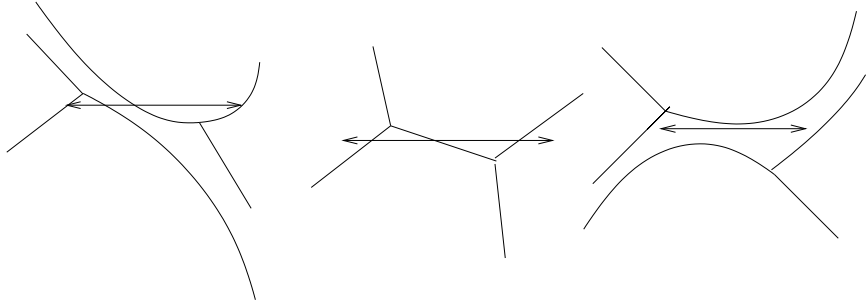}
\end{center}
\caption{Stokes line for $E$ near  $\lambda_0 e^{\frac {i\pi} 4}$, the arrows locate the positions of $-1$ and $1$ }
\label{evo}
\end{figure}

The second statement is proved by  the relation 
$$ S_{0,\alpha_+}(E) +  S_{\alpha_+,1}(E) + S_{1,0}(E)=0,$$\par
and illustrated in Figure \ref{evo}.\par

Point 3 is proved in Theorem \ref{alinfini} below.\par
\begin{figure}[ht]
\begin{center}
\epsfysize=4cm
\epsfbox{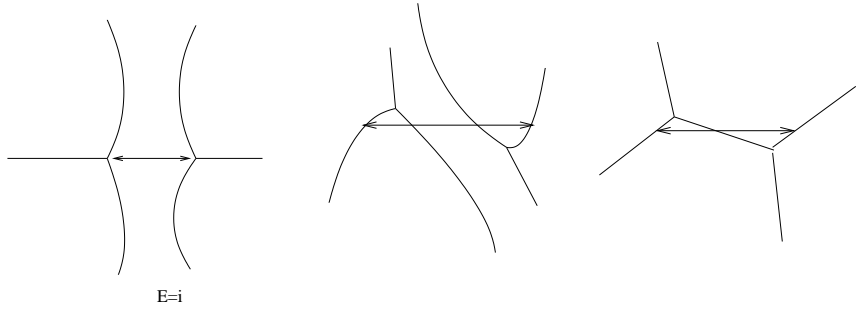}
\end{center}
\caption{Evolution of the Stokes line for $E\in \Gamma_{\alpha_+,1}$, the arrows locate the position of $-1$ and $1$}
\label{extre}
\end{figure}

Point 4 is deduced from
 $$z(x,y,E)=\pm \int_{xe^{\frac{i\pi} 8}} ^{y e^{\frac{i\pi} 8}} (u^2-Ee^{-\frac{i\pi} 4} )^{1/2} \,du$$
 So for $E= \{\lambda e^{\frac {i\pi} 4},\lambda\in\rit\}$ we obtain
  $$\re S_{0,\alpha_+}(E)=\re \pm \int_{0} ^{\pm \sqrt\lambda} (u^2-\lambda )^{1/2} \,du=0.$$
The existence of $\lambda_0$ is obvious in Figure \ref{extre}.\par

Point 5 is obvious from Figure \ref{extre}.\par

\cqfd\par

The hypothesis  \textbf{(H1)} is obviously satisfied by the potential $V(x)=ix^2$.
The hypothesis \textbf{(H2) (H3)} are satisfied by the potential $ix^2$ ( and by any  potential $e^d x^2$  if $e^{\frac d 2}\notin \rit $ ) but not by the potential $x^2$.  
 
For large energies, we summarize the result of Proposition \ref{inver} , Theorem \ref{reci}   and Theorem \ref{curve} in
\begin{theorem}\sl We have\par\noindent
\begin{enumerate}
\item Large values of $\lim_h\Sp P^h$  are  close to the curve $\Gamma_{-1,1}$. 
\item The curve  $\Gamma_{-1,1}=\{E\in \C, \re S_{-1,1}(E)=0\}$ tends  to infinity, and  is asymptotic to the line $\R+i/3$.
\end{enumerate}
\label{alinfini}
\end{theorem}
 One could compute the set $\CT^c$ with Matlab or Mathematica. In Figure \ref{spectre}  we have  drawn the set $\CT^c$ using Theorems \ref{forme} and \ref{alinfini}. \par
\par\sl
For large energies and perturbation , we summarize the results of Theorem \ref{reciper} and Theorem \ref{curveper} in
\begin{theorem}\sl We have\par\noindent
\begin{enumerate}
\item Large  $E\in\lim_h\Sp P_\delta^h$  are  near to   the union of the two  curves
$\Gamma_{-1,\beta}$, $\Gamma_{\beta,1}$. 
\item The curve  $\Gamma_{-1,\beta}=\{E\in \C, \re S_{-1,\beta}(E)=0\}$ tends to infinity, and is asymptotic to  $\R+i \frac 1 3 (\beta^2-\beta+1) -i\delta$.
\item The curve  $\Gamma_{\beta,1}=\{E\in \C, \re S_{\beta,1}(E)=0\}$ goes to infinity, and  is asymptotic to $\R+i \frac 1 3 (\beta^2+\beta+1) +i\delta$.
\end{enumerate}
\label{infiniper}
\end{theorem}
We remark that the two curves $\Gamma_{-1,\beta}$ $\Gamma_{\beta,1}$ are distinct for $\beta\neq0$.

Figure \ref{spectreper} represents the two Theorems,   \ref{alinfini} and \ref{infiniper}.\par

\par\sl

%

\section{pseudo-spectrum}
\label{speudospec}
We have defined here the pseudo-spectrum of $H_0^h$ as the set 
$${\overline{ \{\xi^2+V(x), x\in(-1,1), \xi\in\rit, \im(\xi V'(x))\neq 0\}}}.$$ 
We remark that for any $z$ is this set, we could contruct a function $u^h$ in the domain of $H_0^h$  such that 
$$\|(H^h_0-z)u^h\|=\O(h^\infty) \|u^h\|.$$

Let suppose now that we put two jumps in the potential, one at  $\beta$ and one at $\beta'$. 
i.e $V_{\beta\beta'}=V+i\delta_1 H_\beta+i\delta_2 H_{\beta'}.$ 

Assume that $V$ verifies the  hypothesis  \textbf{(H1)} and~:  

\textbf{(H4)} If  $E$ is  large enough, with small enough imaginary part then  $1$,  $-1$, $\beta$ and $\beta'$  belong to the same Stokes region.

Then we obtain just as in Theorems \ref{reciper} and \ref{curveper}  the next two theorems 
\begin{theorem}
For $E_0$ large, the following two conditions are equivalent
\begin{enumerate}
\item $\re(\int_{-1}^\beta\sqrt {V_{\beta\beta'}-E_0}\,dx=0)$ or $\re(\int_{\beta'}^1\sqrt {V_{\beta\beta'}-E_0}\,dx=0)$ or \par\noindent
$\re(\int_{\beta}^{\beta'}\sqrt {V_{\beta\beta'}-E_0}\,dx=0)$ .
\item $E_0\in\lim_h\Sp H^h_\delta$.
\end{enumerate}
\end{theorem}

\begin{theorem}
For  $a$ large enough, the equation
$$
\re S_{\beta,\beta'}(E)|_{E=a+ib}= \re \int_{\beta}^{\beta'} \sqrt{V(x)+i\delta_1-a-ib}\,dx=0
$$
has  a unique solution $b(a)$, the solution  satisfies  
$$ b(a)=i \frac 1 {\beta'-\beta} \im(Y(\beta')-Y(\beta))-i \delta_1 +O( \frac 1 a ).$$
\end{theorem}
We remark that 
$$\lim_{\beta'\to \beta}\frac 1 {\beta'-\beta} \im(Y(\beta')-Y(\beta))=\im(V(\beta))$$
This means, letting $\beta'\to\beta$, and $\delta\to 0$ any values of the form 
$\{\xi^2+i\im( V(x)); \xi \in\rit, |\xi|>>1, x\in[-1,1] \}$ is  in the
spectrum of this kind of  perturbation of $H^h$.
 Thus we can reach all the large values of the  pseudo-spectrum  with this special type of perturbation.

%

I wish to thank J. Sj\"ostrand  for  organizing a working group  in IHP,
where  one of the subjects was  pseudo-spectrum. I would also like to thank
T. Ramond, who motivated my interest in this problem and helped me to obtain this
final presentation, T. Kappeler for useful discussions about
the history of Sturm Liouville problems, M. Zworski for pointing out  the paper
\cite{mi}. Finally, I wish to thank the ETH in Z\"urich where
part of this work was done.\end{document}